\documentclass[11pt]{amsart}

\let\oldtocsection=\tocsection

\let\oldtocsubsection=\tocsubsection

\let\oldtocsubsubsection=\tocsubsubsection

\renewcommand{\tocsection}[2]{\hspace{0em}\oldtocsection{#1}{#2}}
\renewcommand{\tocsubsection}[2]{\hspace{1em}\oldtocsubsection{#1}{#2}}
\renewcommand{\tocsubsubsection}[2]{\hspace{2em}\oldtocsubsubsection{#1}{#2}}

\usepackage{geometry}                
\usepackage{graphicx}
\usepackage{tikz}
\usetikzlibrary{matrix,arrows,decorations.pathmorphing}
\usepackage{amssymb}
\usepackage{epstopdf}
\title{Deligne Pairings and Discriminants of Algebraic Varieties }

\author{H. Manilal Kapadia}

\begin{document}

\begin{abstract}
Let $V$ be a finite dimensional complex vector space, $V^*$ its dual and let $X\subset\mathbb{P}(V)$ be a smooth projective variety of dimension $n$ and degree $d\geq 2$. For a generic $n$-tuple of hyperplanes $(H_1,...,H_n)\in \mathbb{P}(V^*)^n$, the intersection 
$X\cap H_1\cap \cdots \cap H_n$ consists of $d$  distinct points. We define the ``discriminant of $X$", to be  the the set $D_X$ of $n$-tuples for which the set-theoretic intersection is not equal to $d$ points. Then $D_X\subset \mathbb{P}(V^*)^n$ is a hypersurface and the set of defining
polynomials, which is a one-dimensional vector space, is called the ``discriminant line". We show that this line is
canonically isomorphic to the Deligne pairing $\langle KL^n,...,L\rangle$ where $K$ is the canonical line bundle of $X$ an
$L\rightarrow X$ is the restriction of the hyperplane bundle. As a corollary, we obtain a generalization of Paul's formula
\cite{Paul} which relates the Mabuchi K-energy on the space of Bergman metrics to $\Delta_X$, the ``hyperdiscriminant of $X$". 

\end{abstract}
\newcommand*{\DashedArrow}[1][]{\mathbin{\tikz [baseline=-0.25ex,-latex, dashed,#1] \draw [#1] (0pt,0.5ex) -- (1.3em,0.5ex);}}%

\maketitle

\renewcommand{\thefootnote}{\fnsymbol{footnote}}
\newcommand{\starttext}{ \setcounter{footnote}{0}
\renewcommand{\thefootnote}{\arabic{footnote}}}
\renewcommand{\theequation}{\thesection.\arabic{equation}}
\newcommand{\be}{\begin{equation}}
\newcommand{\bea}{\begin{eqnarray}}
\newcommand{\eea}{\end{eqnarray}} \newcommand{\ee}{\end{equation}}
\newcommand{\N}{{\cal N}} \newcommand{\<}{\langle}
\renewcommand{\>}{\rangle}
\def\ba{\begin{eqnarray}}
\def\ea{\end{eqnarray}}
\newcommand{\PSbox}[3]{\mbox{\rule{0in}{#3}\includegraphics{#1}
\hspace{#2}}}

\def\v{\vskip .1in}

\def\al{\alpha}
\def\b{\beta}
\def\c{\chi}
\def\d{\delta}
\def\e{\epsilon}
\def\g{\gamma}
\def\l{\lambda}
\def\m{\mu}
\def\n{\nu}
\def\o{\omega}
\def\f{\phi}
\def\r{\rho}
\def\si{\sigma}
\def\t{\theta}
\def\z{\zeta}

\def\G{\Gamma}
\def\D{\Delta}
\def\O{\Omega}
\def\T{\Theta}

\def\cA{{\mathcal A}}
\def\cB{{\mathcal B}}
\def\cC{{\mathcal C}}
\def\cD{{\mathcal D}}
\def\cE{{\mathcal E}}
\def\cF{{\mathcal F}}
\def\cG{{\mathcal G}}
\def\cH{{\mathcal{H}}}
\def\cI{{\mathcal I}}
\def\cK{{\mathcal K}}
\def\cL{{\mathcal L}}
\def\cM{{\mathcal M}}
\def\cN{{\mathcal N}}
\def\cO{{\mathcal O}}
\def\cP{{\mathcal P}}
\def\cR{{\mathcal R}}
\def\cS{{\mathcal S}}
\def\cT{{\mathcal T}}
\def\cX{{\mathcal X}}

\def\N{\mathbb N}
\def\Z{{\mathbb Z}}
\def\Q{{\mathbb Q}}
\def\R{{\mathbb R}}
\def\C{{\mathbb C}}
\def\P{{\mathbb P}}
\def\K{{\rm K\"ahler }}

\def\KE{{\rm K\"ahler-Einstein }}
\def\KEE{{\rm K\"ahler-Einstein}}
\def\Ric{{\rm Ric}}
\def\Hom{{\rm Hom}}
\def\mod{{\ \rm mod\ }}
\def\Aut{{\rm Aut}}
\def\End{{\rm End}}
\def\osc{{\rm osc}}
\def\vol{{\rm vol}}
\def\Vol{{\rm Vol}}
\def\ti\tilde

\def\u{\underline}

\def\pl{\partial}
\def\na{\nabla}
\def\i{\infty}
\def\I{\int}
\def\p{\prod}
\def\s{\sum}
\def\dd{{\bf d}}
\def\ddb{\partial\bar\partial}
\def\sub{\subseteq}
\def\ra{\rightarrow}
\def\hra{\hookrightarrow}
\def\Lra{\Longrightarrow}
\def\lra{\longrightarrow}
\def\LA{\langle}
\def\RA{\rangle}
\def\L{\Lambda}

\def\us{{\underline s}}
\def\Re{{\rm Re}}
\def\Im{{\rm Im}}
\def\tr{{\rm tr}}
\def\det{{\rm det}}
\def\half{ {1\over 2}}
\def\third{{1 \over 3}}
\def\ti{\tilde}
\def\un{\underline}
\def\Tr{{\rm Tr}}

\def\pz{\partial _z}
\def\pv{\partial _v}
\def\pw{\partial _w}
\def\w{{\bf w}}
\def\x{{\bf x}}
\def\y{{\bf y}}
\def\z{{\bf z}}
\def\tet{\vartheta}
\def\dwplus{\D _+ ^\w}
\def\dxplus{\D _+ ^\x}
\def\dzplus{\D _+ ^\z}
\def\chiz{{\chi _{\bar z} ^+}}
\def\chiw{{\chi _{\bar w} ^+}}
\def\chiu{{\chi _{\bar u} ^+}}
\def\chiv{{\chi _{\bar v} ^+}}
\def\os{\omega ^*}
\def\ps{{p_*}}

\def\hO{\hat\Omega}
\def\ho{\hat\omega}
\def\o{\omega}

\def\[{{\bf [}}
\def\]{{\bf ]}}
\def\Rd{{\bf R}^d}
\def\Ci{{\bf C}^{\infty}}
\def\pl{\partial}
\newcommand{\dotcup}{\ensuremath{\mathaccent\cdot\cup}}

\parindent=0in

\newtheorem{theorem}{Theorem}
\newtheorem{corollary}{Corollary}
\newtheorem{lemma}{Lemma}
\newtheorem{definition}{Definition}
\newtheorem{proposition}{Proposition}

\setcounter{equation}{0}

\section{Introduction}

Let $V$ be a finite dimensional complex vector space, $V^*$ its dual, and $X\sub\P(V)$ be a smooth projective manifold of dimension $n$ and degree $d\geq 2$. 
\begin{definition}\label{disc}
The discriminant variety of $X$ is the projective variety $D(X)$ defined by
\be\label{D}
D(X)\ = \ \{(H_1,...,H_n)\in \P(V^*)^n: \#\{X\cap H_1\cdots\cap H_n\}\not=d\}
\ee
\end{definition}
Thus, $(H_1,...,H_n)\in D(X)$ if and only if there exists a point $x\in X\cap H_1\cdots\cap H_n$ such that
$\dim(H_1\cap\cdots H_n\cap ET_xX)\geq 1$, where $ET_xX\sub\P(V)$ is the imbedded tangent space
$x$. 
\v
An essential property of $D(X)\sub \P(V^*)^n$ is that it has codimension one. Thus we may consider the one-dimensional vector space $\hat D(X)$ of polynomials of degree $(d,d,...,d)$ which vanish
on $D_X$. We call $\hat D(X)$  the discriminant line.
\v
If $\pi:\cX\ra S$ is a flat family of varieties of dimension $n$ and degree $d$ over a base $S$, then one
can similarly define $\hat\cD_X\ra S$, a line bundle whose fiber at $s\in S$ equals the one dimensional
vector space $\hat D(X_s)$ where $X_s=\pi^{-1}(x)$.
\v
It is useful to compare the definition of $D(X)$ to that of $C(X)$, the Chow variety:
\begin{definition} The Chow variety of $X$ is the projective variety $C(X)$ defined as follows
\be\label{C}
C(X)\ = \ \{(H_0,...,H_n)\in \P(V^*)^{n+1}: \#\{X\cap H_0\cdots\cap H_n\}\not=0\}
\ee
\end{definition}
The variety $C(X)\sub\P(V^*)^{n+1}$ also has codimension one and the one-dimensional vector space of polynomials
of degree $(d,d,...,d)$ which vanish on $C(X)$ is called the Chow line.
\v
In the terminology of the work by Gelfand-Kapranov-Zelevinsky \cite{G}, the variety $D_X$ is essentially the ``first higher associated hypersurface of $X$".
\v

Let $K\ra X$ the canonical line bundle, and $L\ra X$ the restriction
of the hyperplane line bundle.  Zhang's theorem \cite{Z} says that the 
Deligne pairing $\<L,....,L\>$ of $n+1$ copies of $L$ is canonically
isomorphic to the the Chow line. 
\v
 In this paper
we show that the
Deligne pairing $\<KL^n,L,L,...L\>$ of $KL^n$  with $n$ copies of $L$
is canonically isomorphic to the discriminant line. As an application,
we conclude, via the theorem of Phong-Ross-Sturm \cite{PRS}, that
the sub-dominant term in the Mumford-Knudsen expansion is equal to
the discriminant of $X$. A second application is to the asymptotics
of the K-energy on the space of Bergman metrics: The Mabuchi bundle
$\cM$, which was defined and studied by Phong-Sturm \cite{Stab, Fut, PS}.
They showed its associated Deligne metric is precisely the Mabuchi
K-energy. Combining this with our result, we show that the K-energy
on the space of Bergman metrics equals the log of the Deligne norm of the Discriminant Point
minus the log of the Deligne norm of the Chow Point (see Corollary 2 for the statement).
This generalizes, and makes more precise, the theorem of Paul \cite{Paul} (in which
the same formula is proved under a restrictive hypothesis).
\section{Deligne Pairings}

\newcommand*{\DDashedArrow}[1][]{\mathbin{\tikz [baseline=-0.25ex,-latex, dashed,#1] \draw [#1] (0pt,0.5ex) -- (1.3em,0.5ex);}}%



We outline the basic theory, following closely the paper of Zhang \cite{Z}.
Let $\pi:\cX\ra S$ be a flat projective morphism of integral schemes defined over $\C$, of pure relative dimension~$n$.
Let $\cL_0,...,\cL_n$ be line bundles on $X$. Then the Deligne pairing 
$\<\cL_0,...,\cL_n\>(\cX/S)\ra S$ is a line
bundle on $S$ which is defined as follows: A section of $\<\cL_0,...,\cL_n\>$  
over a small open set $U\sub S$ is a symbol $\<l_0,...,l_n\>$
where the $l_j:U \DDashedArrow[->,densely dashed    ] \cX $
are rational sections whose divisors have empty intersection. The relation between the symbols is
given as follows: If $f$ is a generic rational function on $U$ and if 
$$ \prod_{i\not=j}{\rm div}(l_i)\ = \ \s_k n_kY_k
$$
is flat over $S$, then 
$$ \<l_0,...,fl_j,...,l_n\>\ = \prod_k{\rm Norm}_{Y_k/S}(f)^{n_k}\,\<l_0,...,l_n\>
$$
We summarize below some of the properties of Deligne pairings which will be needed.\subsection{Projection Formulas.}
\subsubsection{With $n$ pullbacks}\
\v
Let $\phi: X \ra Y$ and $\pi: Y \ra S$. Let $m = \dim(X/Y)$
and $n = \dim(Y/S)$.
\v
Let $\cK_0, ..., \cK_m$ be line bundles on $X$.
Let $\cL_1,..., \cL_n$ be line bundles on $Y$.

\v
Then
\be\label{n}   \<\cK_0, ..., \cK_m, \phi^*\cL_1, ... ,\phi^*\cL_n\>
(X/S)
\ = \
\<\<\cK_0, ..., \cK_m\>, \cL_1, ... ,\cL_n\>(Y/S)
\ee

The map is given by
$F: \<k_0,...,k_m,\phi^*l_1,..., \phi^*l_n\> \ \mapsto\
 \<\<k_0,...,k_m\>,l_1,..., l_n\> $. 
 \subsubsection{With $n+1$ pullbacks}\
\v
Let $\cK_1, ..., \cK_m$ be line bundles on $X$.
Let $\cL_0, \cL_1,..., \cL_n$ be line bundles on $Y$.
\v
Then
\be\label{n1}  \<\cK_1, ..., \cK_m, \phi^*\cL_0, ... ,\phi^*\cL_n\>
(X/S)
\ = \
\<\cL_0, \cL_1, ... ,\cL_n\>(Y/S)^D
\ee
where $D = deg[c_1(\cK_{1\eta})\cdot c_1(\cK_{2\eta})\cdots
c_1(\cK_{m\eta})]$ and $\eta$ is a generic point
on $Y$. In other words, $D$ is the number of points in
$div(\cK_1) \cap \cdots \cap div(\cK_m)$ in a generic fiber.
\v
\subsubsection{With $n+2$ pullbacks}\
\v
Let $\cK_1, ..., \cK_{m-1}$ be line bundles on $X$.
Let $\cL_0, \cL_1,..., \cL_{n+1}$ be line bundles on $Y$.
\v
Then
\be\label{n2}   \<\cK_1, ..., \cK_{m-1}, \phi^*\cL_0, ...
,\phi^*\cL_{n+1}\>  (X/S)
\ = \
\cO_S
\ee
\subsection{Induction Formula}\
\v
Let $\pi: X \ra S$ and $\cL_0, ..., \cL_n$ as before.
Let $l$ be a rational section of $\cL_n$. Assume
all components of $div(l)$ are flat over $S$.
Then
\be\label{ind}  \<\cL_0,...,\cL_n\>(X/S) \ = \
\<\cL_0,...,\cL_{n-1}\>({\rm div}(l)/S)
\ee
\section{Metrics on Deligne Pairing}

Assume that $h_j$ is a smooth metric on $\cL_j$. If $\cX$ is not smooth, then this means 
that there is a smooth manifold $\cX'$, a smooth line bundle $\cL_j'\ra \cX'$, and a
smooth metric $h_j'$ on $\cL_j$ whose restriction to $\cX$ equals $h_j$. Then the
Deligne metric $\<h_0,...,h_n\>$ is a metric on the line bundle $\<\cL_0,...,\cL_n\>$
which is defined inductively by the formula

\be\label{metric} \log\|\<l_0,...,l_n\>\|\ = \ \log\|\<l_0,...,l_{n-1}\>\|\ + \ \I_{\cX/S}\log|l_n|\,\o_0\wedge\cdots\wedge \o_{n-1}
\ee
where $\o_j=-{i\over 2\pi}\ddb\log|l_j|^2$.
\v
The inductive formula (\ref{metric}) implies that (\ref{ind}) is an isometry
\be \<h_0,...,h_n\>\ = \ \<h_0,...,h_{n-1}\>\exp\left(
-\I_{\cX/S}\log|l_n|\,\o_0\wedge\cdots\wedge \o_{n-1}
\right)
\ee If $\phi_0,...,\phi_n$ are smooth functions on $\cX$, formula (\ref{metric}) also implies
\be \<h_0e^{-\phi_0},...,h_ne^{-\phi_n}\>\ = \ \<h_0,...,h_n\>\exp(-E(\phi_0,...,\phi_n))
\ee
where
$$ E(\phi_0,...,\phi_n)\ = \ \s_{j=0}^n\I_{\cX/S} \phi_j\, (\bigwedge_{k<j} \o_{\phi_k})\wedge(\bigwedge_{k>j} \o_k)
$$
and $\o_{\phi_k}=\o_{k}+\small{{i\over 2\pi}} \ddb\phi_k$.
In particular, if $\cL_0=\cdots =\cL_n$ and $h_0=\cdots =h_n$, then, setting $E_\chi(\phi)=E(\phi,...,\phi)$ we obtain
\be\label{AY} E_{\cX}(\phi)\ = \ \s_{j=0}^n\I_{\cX/S}\phi_j \,\o_\phi^j\,\o^{n-j}
\ee
which coincides with the well known Aubin-Yau functional.
\v
A simple consequence of these formulas which will  later be useful is the following:
\v
\begin{proposition}\label{prop1} 
Let $X$ be a smooth projective variety of dimension $n$, $L_0,...,L_n\ra X$  holomorphic line bundles, and $h_j$ a hermitian metric on $L_j$.
Let $G$ be a semi-simple Lie group acting on $L\ra X$, and define  $\phi_j^\si$ by the following formula.
$$ \si^*h_j\ = \ he^{-\phi_j^\si}
$$ 
Then $E(\phi_0^\si,...,\phi_n^\si)=0$. 
\end{proposition}
{\it Proof.} Define $\r: G\ra \C^\times$ by the formula
$$ \<\si^*s_0,...,\si^*s_n\>\ = \ \r(\si)\<\si_0,...,\si_n\>
$$ 
Then 
$\|\<\si^*s_0,...,\si^*s_n\>\|\ = \ |\r(\si)|\cdot\|\<\si_0,...,\si_n\>\|$ so
$E(\phi_0^\si,...,\phi_n^\si)=-\log|\r(\si)|$. On the other hand, $\r:G\ra \C^\times$ is a homomorphism which
must be trivial since $G$ is semi-simple.
\v
Example. Let $X=\P^N$, $L_j=O(1)$ and  $h=h_{FS}$. Then for every $\si\in SL(N+1)$ we have
\be E_{\P^N}(\phi^\si)\ = \ 0
\ee
\section{The Mabuchi Line Bundle}
Let $\pi:\cX\ra S$ be as above, and assume $\cK_{\cX/S}$, the relative canonical bundle, is well defined and
let $h$ be a positively curved metric on $\cL$ with curvature $\o$. Define $h^{-1}_\cK$, which as a metric
on $\cK^{-1}$, by the formula  $h^{-1}_\cK=\o^n$.
Phong-Sturm \cite{Stab, Fut} introduced the Mabuchi line bundle\footnote{Zhang had also deduced this bundle earlier in a 1993 letter to Deligne.}, which is the hermitian  bundle
\be\label{mab} \cM_h\ = \ \<\cK,\cL,...,\cL\>^{1\over c_1(L)^n} \<\cL,...,\cL\>^{-\m\over c_1(L)^n}
\ee
where $c_1(L)^n$ is computed on a generic fiber, and $\m\in\Q$ is uniquely determined by requiring that the metric is scale
invariant, that is, invariant under $h\mapsto \l h$ for $\l$ a positive real number. If
follows from the definitions that $nc_1(K)c_1(L)^{n-1}-\m(n+1)c_1(L)^n=0$ so
$$ \m\ = \ {n\over n+1}{c_1(K)c_1(L)^{n-1}\over c_1(L)^n}
$$

\v
Then, as shown in \cite{Fut},
\be\label{mab} \cM_{he^{-\phi}}\ = \ \cM_h \exp(-\n(\phi_s))
\ee
where for $s\in X$ with $X_s$ a smooth fiber,  $\n(s,\phi)$ is the Mabuchi K-energy $\phi_s=\phi|_{X_s}$.
\v
For our purposes, it is more convenient to rewrite as follows:

\be \label{mab1} \cM_h^{[c_1(L)^n]}\ = \ \<\cK\cL^n,\cL,...,\cL\> \<\cL,...,\cL\>^{-\m-n}
\ee

The theorem of S. Zhang shows that $\<\cL,...,\cL\>$ is canonically isomorphic to the Chow bundle.
We shall use Zhang's approach to prove $\<\cK\cL^n,\cL,...,\cL\>$ is canonically isomorphic to the
Discriminant bundle.

\section{Tangent bundle for projective space}

Let $V$ be a complex vector space of dimension $N+1$ and $\P(V)=\{x\sub V: \dim(x)=1\}$. Let
$O(1)\ra\P(V)$ be the hyperplane line bundle.
\v
 If
$0\not=\hat x\in V$ and $x = \C\cdot \hat x\in\P(V)$ then we have a canonical map
$$ I_{\hat x}: V=T_{\hat x}V\ra T_x\P(V)
$$
Since $\ker(I_{\hat x}) = x$ we see $I_{\hat x}: V/x \ra T_x\P(V)$ is an isomorphism, or equivalently, the map $I_{\hat x}^*: (V/x)^*\ra T_x^*\P(V)$ is an isomorphism.
\v
Let $J_{\hat x}: T_x\P(V)\ra V/x$ be the inverse. Then for $\al\in\C^\times$ we have
$J_{\al\hat x}=\al J_{\hat x}$ so we see $J: O(1)\ra (V/x)^*\otimes T_x\P(V)$ is a vector bundle map
and $O(1)\otimes V/x \ra T_x\P(V)$ is a canonical isomorphism:
\be O(-1)\otimes T_x\P(V)=V/x
\ee
Alternatively,
$$ T_x^*\P(V)\ = \ (V/x)^*\otimes O_x(-1)
$$
\section{The Deligne metric}
Let $Y$ be a projective manifold of dimension $m$, $L\ra Y$ an ample line bundle, and $d$ a positive integer. Let $V=c_1(L)^m$ and
$\cN\ra \P(H^0(Y,L^d))$ be the hyperplane line bundle. 
\v
Let $h$ be
a positively curved hermitian metric on $L$ and $\o=-i\ddb\log h>0$ its curvature. Let ${\rm D}(h)$
be the  norm on the vector space $H^0(Y,L^d)$ defined by the following formula. If $0\not=f\in H^0(Y,L^d)$
then
\be\label{del} \log\|f\|_{D(h)}^2\ = \ {1\over V}\I_Y \log|f|^2_{h^d}\,\o^{m}
\ee
where $V=\I_X 1\,\o^n$.
In particular, ${\rm D}(h)$ makes $\cN$ into a hermitian line bundle.
\v
Remark: Note that ${\rm D}(h)$ is not equal to ${\rm Hilb}(h)$, but 
${\rm D}(h)=e^\psi {\rm Hilb}(h)$ for some bounded smooth function $\psi$ on $H^0(Y,L^d)$ with the property:
$\psi(\l v)=\psi(v)$ for all positive real numbers $\l$.
\v
Let $0\not= f\in H^0(X,L^d)$ and assume $Z=\{f=0\}\sub Y$ is smooth. 
\v
 The map
$$ I_f:\langle L,...,L\rangle_Z \ \ra \ \langle L,...,L^d\rangle_Y
$$
given by $\langle s_0,...,s_{m-1}\rangle_Z\mapsto \langle s_0,...,s_{m -1},f\rangle_Y$ is
an isomorphism. Since $I_{\al f}=\al^{c} I_f$ we see

$$ \cN^{-V}_{[f]} = \langle L,...,L^d\rangle_Y\otimes \langle L,...,L\rangle_Z^{-1}
$$
Equivalently there is a canonical isomorphism
\be\label{gen} J_{[f]}:\cN_{[f]}^V\ \ra \  \langle L,...,L\rangle_Z\otimes \langle L,...,L\rangle_Y^{-d}
\ee
which is easily seen to be an isometry (with metric $D(h)$ on the left and the Deligne
metric on the right). Moreover, $J$ is $G\sub GL(H^0(Y,L^d))$ equivariant, where $G=Aut(Y,L)$.
\v
Let $L\ra Y$ be a holomorphic line bundle on a projective manifold $Y$ and $h$
a smooth metric on $L$.
Suppose $G\sub \Aut(Y,L)$ is a semi-simple Lie group, and write $\si^*h=he^{-\phi_\si}$ for $\si\in G$.
\begin{corollary}\label{cor1} Let $f\in H^0(Y,L^d)$ be such that $Z=\{f=0\}$ is a smooth sub manifold. If $E_Z$ is
the Aubin-Yau functional on $Z$ (\ref{AY}) then for all $\si\in G$ we have
\be E_{Z}(\phi^\si) \ = \ {1\over V}\log\left({\|f^\si\|^2_{D(h)}\over \|f\|^2_{D(h)}}\right)
\ee
\end{corollary}

\section{The theorem of Zhang}
In this section we give a slightly modified version of Zhang's proof.
Let $X^{(n)}\sub \P(V)$  of degree $d$, and write $L\ra \P(V)$ and $M\ra \P(V^*)$ for the hyperplane bundles.
Let $\P=\P(V^*)^{n+1}$ and $M_i=\pi_i^*M\ra\P$. Let
$$\cM=\pi_1^*M\otimes\cdots\otimes \pi_{n+1}^*M=M_1\otimes\cdots\otimes M_{n+1}\ra\P$$
Let $\pi_X: X\times\P\ra X$ and $\pi_\P:X\times\P\ra\P$ be the projection maps and consider
\be
\label{big}
\cB=\langle 
\overbrace{
\hbox{$ \pi_\P^*\cM,,...\pi_\P^*\cM$}
}^
{
\hbox{$m$}
},
\overbrace{
\hbox{$
 \,\pi_X^*L\otimes\pi_\P^*M_1,...,
\pi_X^*L\otimes \pi_\P^*M_{n+1}$}
}^
{
\hbox{$n+1$}
}
\rangle_{X\times\P/*} \ee
where $m=\dim\P$ and $*$ is a point.
We evaluate $\cB$ in two different ways. First,
let 
$$ \G\ = \ \{(x,H_1,....,H_{n+1})\in X\times\P\ : \ x\in H_1\cap\cdots\cap H_{n+1}\ \}
$$
and let $Z=\pi_\P(\G)$. Then $Z=C(X)\sub \P$ is the Chow hypersurface of~$X$. 
The line $\langle \cM,...,\cM\rangle_{Z/*}=\langle \cM,...,\cM\rangle_{Z}$ is called the Chow line.
\v
Next we define a section $s_i$ of $\pi_X^*L\otimes\pi_\P^*M_i$ as follows: 
$s_i(x,H) \in x^*\otimes H^* = \Hom(x\otimes H,\C)$ is the restriction of the canonical
paring $\Hom(V\otimes V^*,\C)$, in other words
\be\label{s} s_i(x,H)(z,\l)\ = \ \l(z)
\ee
for all $z\in x$ and $\l\in H$. Note that $s_i(x,H)=0$ if and only if $x\in H$. Thus, applying
(\ref{ind}) a total of $n+1$ times we obtain:

\be\label{first}
\cB\ = \ \langle 
\overbrace{
\hbox{$\cM,\cM,...,\cM$}
}^
{\hbox{$m$}}
\rangle_{Z/*}
\ee

On the other hand, expanding the last $n+1$ terms on (\ref{big}) we have $\cB=\cB_1\otimes\cB_2$
where
\be
\label{big1}
\cB_1=\langle 
\overbrace{
\hbox{$ \pi_\P^*\cM,...\pi_\P^*\cM$}
}^
{
\hbox{$m$}
},
\overbrace{
\hbox{$
 \,\pi_X^*L,...,
\pi_X^*L$}
}^
{
\hbox{$n+1$}
}\rangle_{X\times\P}
\ee

and

\be
\label{big2}
\cB_2=\prod_{i=1}^{n+1}\langle 
\overbrace{
\hbox{$ \pi_\P^*\cM,...,\pi_\P^*\cM$}
}^
{
\hbox{$m$}
},
\overbrace{
\hbox{$
 \,\pi_X^*L\otimes\pi_\P^*M_1,...,
\pi_\P^*M_i,...
\pi_X^*L\otimes \pi_\P^*M_{n+1}$}
}^
{
\hbox{$n+1$}
}\rangle_{X\times\P}
\ee
Now (\ref{n1}) gives
\be\label{B1}
\cB_1\ = \ \langle L,...,L\rangle_{X}^{\deg(\cM)}
\ee
If we expand the last $n+1$ terms in $\cB_2$, using (\ref{n2})the only terms which survive are
those of the form

$$
\overbrace{
\hbox{$
\langle\pi_\P^*M_i, \pi_X^*L,
...
\pi_X^*L$}
}^
{
\hbox{$n+1$}
}\rangle_{X\times\P}
$$
since, by (\ref{n2}), the other terms vanish. Applying (\ref{n1}) once again:

\be
\label{big2}
\cB_2=\prod_{i=1}^{n+1}\langle 
\overbrace{
\hbox{$ \cM,...,\cM$}
}^
{
\hbox{$m$}
}
M_i
\rangle_{\P}^{d}\ = \ 
\overbrace{
\hbox{$\langle \cM,...,\cM$}
}^
{\hbox{$m+1$}}
\rangle_{\P}^{d}
\ee
where $d=c_1(L)^n$.
Since $\deg(\cM)=1$ we conclude

$$\langle L,...,L\rangle_{X}\ = \ \langle \cM,...,\cM\rangle_{C(X)}\otimes
\langle \cM,...,\cM\rangle_{\P}^{-d}
$$
Combining with (\ref{gen}) we obtain Zhang's theorem.

\v\v

\section{Deligne pairings and discriminants}
Let $V$ be a vector space over $\C$ of dimension $N+1$ and $X\sub\P(V)$ a smooth projective
variety of degree $d>1$. Recall that  $D(X)\sub \P=\P(V^*)^n$ is a hypersurface. 
Let $L\ra \P(V)$ and $M\ra\P(V^*)$ be the hyperplane line bundles, let $M_j=p_j^*M\ra\P$ and  let $\cM\ra\P$ be the line bundle $\cM= M_1\otimes\cdots\otimes M_n$. 
\begin{theorem}
$$ \<KL^n,L,...,L\>_{X} = \<\cM,...,\cM\>_{D(X)}\otimes\langle \cM,\cM,...,\cM\rangle^{-d}_{\P}
$$
\end{theorem}
{\it Proof.}
Let 
$m=\dim(\P)$ and

$$
\cB=\langle 
\overbrace{
\hbox{$ \pi_\P^*\cM,,...\pi_\P^*\cM$}
}^
{
\hbox{$m$}
}, \pi_X^*(K\otimes L^n)\otimes\pi_\P^*\cM, \pi_X^*L\otimes\pi_\P^*M_1,...,\pi_X^*L\otimes\pi_\P^*M_n
\rangle_{X\times\P} \ \longrightarrow \ X\times\P
$$
\newpage
On the one hand,

$$ \cB\ = \ \langle \cM,\cM,...,\cM\rangle^d_{\P}\otimes \<KL^n,L,...,L\>_{X}
$$
where 
$$ d\ = \ (n+1)c_1(L)^n+ c_1(K)c_1(L)^{n-1}
$$
On the other hand,

$$ \cB\ = \ \langle \pi_\P^*\cM,...,\pi_\P^*\cM, \pi_X^*(K\otimes L^n)\otimes\pi_\P^*\cM\rangle_{\G'}
$$
where 
$$\G' \ = \ \{(x,H_1,...,H_n)\in X\times\P: x\in H_1\cap\cdots\cap H_n\}
$$
Now we wish to define a holomorphic section $s: \G'\ra \pi_X^*(K\otimes L^n)\otimes\pi_\P^*\cM$. We define
$$ s(x,H_1,...,H_n)\in \Lambda^n(T_xX^*)\otimes O(n)_x\otimes H_1^*\otimes \cdots H_n^*
$$
$$ \ = \ (\Lambda^n(T_xX)\otimes O_x(-1)\otimes H_1,...,O_x(-1)\otimes H_n)^*
$$
as follows: Since $O_x(-1)\otimes H_j \sub T_X\P(V)^*$, if $\eta_j\in O_x(-1)\otimes H_j$
then $\eta_j\in T_x\P(V)^*$ so $\eta_1\wedge\cdots\wedge\eta_n\in \Lambda^nT_x\P(V)^*$. Hence
if $\o\in\Lambda^nT_xX$ we have a canonical multilinear map
$$ [\Lambda^n(T_xX)]\times  [O_x(-1)\otimes H_1]\times\cdots [O_x(-1)\otimes H_n]\ \ra\C
$$
given by
$$ (\o,\eta_1,...,\eta_n)\mapsto [\eta_1\wedge\cdots\wedge\eta_n](\o)
$$
Thus
$$ \cB\ = \ \<\cM,...,\cM\>_D
$$
where $D=\{s=0\}\sub \G'$. This proves the theorem.
\v
Combining Theorem 1 and (\ref{mab}) we obtain:
\begin{corollary}\label{cor2}
 Let $X\sub \P^N$ be a smooth variety of dimension $n$ and degree $d\geq 2$. Then
\be\label{Kenergy}
\n_\o(\phi_\si)\ = \ \deg(C_X)\log{\|D_X^\si\|^2\over \|D_X\|^2}\ - \ \deg(D_X)\log{\|C_X^\si\|^2\over \|C_X\|^2}
\ee
where $\n_\o$ is the Mabuchi K-energy and the norm is the Deligne norm defined by (\ref{del}).
\end{corollary}
Remark: Suppose that $X\sub\P^N$ is a smooth variety as above whose dual defect vanishes (which holds in the case where $X$ is linearly normal). Then Paul \cite{Paul}
defines $\D_X$, the X-hyperdiscriminant of $X$ to be the dual variety of $X\times \P^{n-1}$, viewed as a sub variety of
$\P^{nN+n-1}$ via the Segre imbedding. Using different methods, he proves  formula (\ref{Kenergy}) with $D_X$ replaced by $\D_X$ and with the Deligne norm replaced by an inexplicit norm. It is
not hard to show that $D_X$ and $\D_X$ are canonically isometric in the case where the dual defect vanishes. Thus
Corollary \ref{cor2} may be viewed as a generalization of Paul's theorem: we don't place any requirement on the dual defect of $X$. Also, our norm is explicitly given by formula (\ref{del}).
\v
Remark: In order to simplify notation, we have restricted to the case where the base $S$ is a single point. One can generalize our results to the case where the base is arbitrary: The Chow line and the Discirminant line, which are one dimensional vector spaces when $S$ is a point, become line bundles on $S$ when $S$ has positive dimension.

\v

{\bf Acknowledgement}. I would like to thank my adviser, Professor Jacob Sturm, for all his help and guidance throughout my years as a graduate student.

\v\v
\bibliographystyle{unsrt}

\vskip .4in
{\it e-mail address:}\ \tt{{hetalk@rutgers.edu}}
\v
\sc{Department of Mathematics,
Rutgers University, 
Newark, NJ 07102}
\end{document}